\documentclass{amsart}
\usepackage{enumerate, amsfonts, amssymb,xypic,amsmath,amsthm}

\newlength{\tabwidth}
\newlength{\tabheight}
\newlength{\tabrule}
\newlength{\tabwidthx}
\newlength{\tabheightx}

\def\gentabbox#1#2#3#4{\vbox to \tabheight{\setlength{\tabrule}{#3}%
  \setlength{\tabwidthx}{#1\tabwidth}\addtolength{\tabwidthx}{\tabrule}%

\setlength{\tabheightx}{#2\tabheight}\addtolength{\tabheightx}{-\tabheight}%
  \hbox to #1\tabwidth{%
 \hspace{-0.5\tabrule}\rule{\tabrule}{#2\tabheight}\hspace{-\tabrule}%
    \vbox to #2\tabheight{\hsize=\tabwidthx%
      \vspace{-0.5\tabrule}\hrule width\tabwidthx height\tabrule%
      \vspace{-0.5\tabrule}\vfil%
      \hbox to \tabwidthx{\hss#4\hss}%
        \vfil\vspace{-0.5\tabrule}%
      \hrule width\tabwidthx height\tabrule\vspace{-0.5\tabrule}}%
 \hspace{-\tabrule}\rule{\tabrule}{#2\tabheight}\hspace{-0.5\tabrule}}%
  \vspace{-\tabheightx}}}
\def\genblankbox#1#2{\vbox to \tabheight{\vfil\hbox to
#1\tabwidth{\hfil}}}
\def\tabbox#1#2#3{\gentabbox{#1}{#2}{0.4pt}{\strut #3}}

\catcode`\:=13 \catcode`\.=13 \catcode`\;=13
\catcode`\>=13 \catcode`\^=13
\def:#1\\{\hbox{$#1$}}
\def.#1{\tabbox{1}{1}{$#1$}}
\def>#1{\tabbox{2}{1}{$#1$}}
\def^#1{\tabbox{1}{2}{$#1$}}
\def;{\genblankbox{1}{1}\relax}
\catcode`\:=12 \catcode`\.=12 \catcode`\;=12
\catcode`\>=12 \catcode`\^=7

\newenvironment{tableau}{\bgroup\catcode`\:=13 \catcode`\.=13
  \catcode`\;=13 \catcode`\>=13 \catcode`\^=13
  \setlength{\tabheight}{3ex}\setlength{\tabwidth}{3ex}%
  \def\b##1##2##3{\gentabbox{##1}{##2}{1.2pt}{\vbox{##3}}}%
  \def\n##1##2##3{\gentabbox{##1}{##2}{0.4pt}{\vbox{##3}}}%
  \vbox\bgroup\offinterlineskip}{\egroup\egroup}

\newtheorem{theorem}{Theorem}[section]
\newtheorem{corollary}[theorem]{Corollary}
\newtheorem{lemma}[theorem]{Lemma}
\newtheorem{proposition}[theorem]{Proposition}

\newtheorem*{theorem*}{}

\theoremstyle{definition}
\newtheorem{definition}[theorem]{Definition}

\theoremstyle{remark}

\newtheorem{example}[theorem]{Example}

\begin{document}
\title{A Relation for Domino Robinson-Schensted Algorithms}
\author{Thomas Pietraho}
\email{tpietrah@bowdoin.edu} 
\address{Department of Mathematics\\Bowdoin College\\Brunswick,
Maine 04011}
\begin{abstract}
We describe a map relating hyperoctahedral Robinson-Schensted
algorithms on standard domino tableaux of unequal rank.  Iteration
of this map relates the algorithms defined by Garfinkle and
Stanton--White and when restricted to involutions, this construction
answers a question posed by M.\! A.\! A.\! van Leeuwen. The principal
technique is derived from operations defined on standard domino
tableaux by D.\! Garfinkle which must be extended to this more general setting.
\end{abstract}
\maketitle
\section{Introduction}

The classical Robinson-Schensted algorithm defines a bijection
between the elements of the symmetric group $S_n$ and same-shape
pairs of standard Young tableaux of size $n$.  The work of Garfinkle
\cite{garfinkle1} defines  similar bijections for $H_n$, the
hyperoctahedral group on $n$ letters, using pairs of certain
same-shape standard {\it domino} tableaux as parameter sets.

Viewing $H_n$  as the Weyl group of a simple Lie group of type $C$,
Garfinkle's generalization is a map $G_0$ whose image is precisely
the set of same-shape pairs of standard domino tableaux of size $n$
and rank $0$. When viewing $H_n$  as the Weyl group of a simple Lie
group of type $B$, she defines a more natural map $G_1$ whose image
is the set of same-shape pairs of standard domino tableaux of size
$n$ and rank one. M.\! A.\! A.\! van Leeuwen has observed that
Garfinkle's definition can be extended to define bijective maps
$G_r$ from $H_n$ to same-shape pairs of standard domino tableaux of
arbitrary rank $r$ \cite{vanleeuwen:rank}.  For $r$ sufficiently
large, $G_r$ recovers the bijection of Stanton and White defined
between $H_n$ and pairs of same-shape standard bitableaux (cf.
\cite{stanton:white} and also \cite{okada}).

Consider an element $\sigma \in H_n$ and let $(T,S) = G_r(\sigma)$
and $(T',S')= G_{r+1}(\sigma)$.  The main result of this paper
describes a map between the pairs $(T,S)$ and $(T',S')$ using
techniques from \cite{garfinkle1}.  In this way, we obtain maps that
relate the different members of this family of generalized Robinson-Schensted algorithms as well as the algorithm and Stanton and White.  When
$\sigma$ is an involution, the map sending $(T,S)$ to $(T',S')$ has
a particularly simple description and answers a question posed by
M.\! A.\! A.\! van Leeuwen in \cite{vanleeuwen:rank}, p. 26.

The combinatorial results of this paper are particularly relevant to recent results in the study of the Kazhdan-Lusztig cell structure of an unequal parameter Hecke algebra $\mathcal{H}$.  Garfinkle's original work on the primitive spectrum of a universal enveloping algebra of a complex semisimple Lie algebra classified the Kazhdan-Lusztig cell structure of equal parameter Hecke algebras of type $B_n$.  In the more general setting of unequal parameter $\mathcal{H}$, \cite{bgil} conjectures a parametrization of cells via domino tableaux of rank $r$, where the specific choice of $r$ depends on the underlying parameters of $\mathcal{H}$. In \cite{pietraho:unequal}, the results of the present paper are used to reconcile the above conjecture and Garfinkle's original work on primitive ideals.  In related work, \cite{gordon:calogero} and \cite{gordon:quiver} provide a geometric interpretation of these combinatorial results in the setting of rational Cherednik algebras.

\section{Definitions and Preliminaries}

\subsection{Generalized Robinson-Schensted Algorithms}

Following Garfinkle \cite{garfinkle1}, we view the elements of the
hyperoctahedral group $H_n$ as those subsets $\sigma$ of
$\mathbb{N}_n \times \mathbb{N}_n \times \{\pm 1\}$, with
$\mathbb{N}_n=\{1, 2, \ldots, , n\}$, such that the projections onto
the first and second components of $\sigma$ are always bijections
onto $\mathbb{N}_n$ (\cite{garfinkle1}, (1.1.2)). We will write the
element $\sigma$ as $\{(\sigma_1,1,\epsilon_1), \ldots,
(\sigma_n,n,\epsilon_n)\}$. In this form, $\sigma$ corresponds to
the signed permutation $(\epsilon_1 \sigma_1, \epsilon_2 \sigma_2,
\ldots, \epsilon_n \sigma_n)$.

For us, Young diagrams will be finite left-justified arrays of
squares arranged with non-increasing row lengths.  A square in row
$i$ and column $j$ of the diagram will be denoted $S_{i,j}$ so that
$S_{1,1}$ is  the uppermost left square in the Young diagram below:

$$
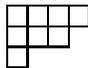
\begin{tiny}
\begin{tableau}
    :.{}.{}.{}.{}\\
    :.{}.{}.{}\\
    :.{}\\
\end{tableau}
\end{tiny}
$$

\begin{definition}
    Let $r \in \mathbb{N}$ and $\lambda$ be a partition of a positive integer $m$.
     A {\it domino tableau of rank $r$ and shape $\lambda$} is a Young diagram
    of shape $\lambda$ whose squares are labeled by integers from some set $M$ in such a way that $0$ labels the square
    $(i,j)$ iff $i+j<r+2$, each element of $M$ labels exactly two adjacent squares, and all labels increase weakly along both rows and columns. A domino tableau is {\it standard} iff $M=\mathbb{N}_n$ for some $n$.
\end{definition}

We will write $DT_r(\lambda)$ for the family of all domino tableaux
of rank $r$ and shape $\lambda$  and $DT_r(n)$ for the family of
all domino tableaux of rank $r$ which contain exactly $n$ dominos.  The corresponding families of standard tableaux will be denoted $SDT_r(\lambda)$ and $SDT_r(n)$.
The set of squares in a tableau $T$ labeled by the integer $l$ will be denoted by $supp(l,T)$ and $supp(0,T)$ will be called the {\it core} of $T$.

Following \cite{garfinkle1} and \cite{vanleeuwen:rank}, we 
describe the Robinson-Schensted bijections
$$G_r : H_n \rightarrow SDT_r(n) \times SDT_r(n) .$$ The algorithm
is based on an insertion map $\alpha$ which, given an element $(i,j,
\epsilon)$ of $\sigma \in H_n$, inserts a domino with label $i$ into
a domino tableau. 

\begin{definition} Consider $\sigma \in H_n$, $(i,j,\epsilon) \in \sigma$, and a domino tableau $T' \in DT_r(k)$.  Write $\ell=\{l_1, l_2,\ldots, l_k\}$ for the set of labels of the dominos of $T'$ listed in increasing order.  When $i \notin \ell$, we can define a tableau $T=\alpha((i,j,\epsilon),T') \in DT_r(k+1)$ by the following procedure:

    \begin{enumerate}
        \item If $i>l_k$, $T$ is formed by:
            \begin{enumerate}
                \item adding a new horizontal domino with label $i$ to the end of the first row of $T'$ if $\epsilon = 1$, or by
                \item adding a new vertical domino with label $i$ at the end of the first column of $T'$ if $\epsilon = -1$.
            \end{enumerate}
        \item Otherwise,  let $l_m$ be the least label in $\ell$ greater than $i$.   We inductively define a sequence $\{T_{m-1}, T_m, \ldots, T_{k+1}\}$ of domino tableaux and let $T=T_{k+1}$. To this effect, construct  $T_{m-1}$ by removing all dominos with labels greater or equal to $l_{m}$ from $T'$.  Let $T_m= \alpha((i,j,\epsilon),T_{m-1})$.
            For $p\geq m$, 
                \begin{enumerate}
                    \item if $supp(l_p,T') \cap T_{p} = \emptyset$, then $T_{p+1}$ is the tableau obtained from $T_{p}$ by labeling $supp(l_p,T')$ with the integer $l_p$;
                    \item if $supp(l_p,T') \cap T_{p} = \{S_{ij}\}$, then
                        $T_{p+1}$ is the tableau obtained from $T_{p}$ by labeling $\{S_{i,j+1},S_{i+1,j+1}\}$ with the integer $l_p$ if $supp(l_p,T')$ is horizontal, or by labeling $\{S_{i+1,j},S_{i+1,j+1}\}$ with the integer $l_p$ if $supp(l_p,T')$ is vertical.
                    \item if $supp(l_p,T') \cap T_{p} = supp(l_p,T')$, then $T_{p+1}$ is the tableau obtained by adding a horizontal domino with label $l_p$ at the end of row $\iota+1$ of $T_{p}$ if $supp(l_p,T')$ is horizontal and lies in row $\iota$ of $T'$, or by adding a vertical domino with label $l_p$ at the end of column $\iota+1$ of $T_{p}$ if $supp(l_p,T')$ is vertical and lies in column $\iota$ of $T'$.
                \end{enumerate}

    \end{enumerate}
\end{definition}

That this procedure is well-defined and indeed produces a domino tableau is verified in \cite{garfinkle1}, Section 2.  To describe the generalized Robinson-Schensted algorithm itself, we start by constructing the left tableau. 
Let $T(0)$ be the only tableaux in $SDT_r(0)$. Define $T(1) = \alpha
((\sigma_1,1,\epsilon_1),T(0))$ and continue inductively by letting
$$T(k+1)=\alpha\big((\sigma_{k+1},k+1,\epsilon_{k+1}),T(k)\big).$$  The left domino tableau
$T(n)$ will be standard and of rank $r$.  The right tableaux track
the shapes of the left tableaux.  Begin by forming a domino tableau
$S(1)$ by adding a domino with label $1$ to $T(0)$ in such a way
that $S(1)$ and $T(1)$ have the same shape. Continue adding dominos
by requiring that at each step $S(k)$ lie in $SDT_r(k)$ and have the
same shape as $T(k)$.  Again, the domino tableau $S(n)$ will be
standard and of rank $r$.  Finally, the image of $\sigma$ under
$G_r$ is defined as the tableau pair $(T(n),S(n))$. To simplify
notation, we will write $G_r^k(\sigma)$ for the pair $(T(k),S(k))$.
We will also sometimes simplify notation slightly and write
$\alpha_{m}(T)$ instead of
$\alpha\big((\sigma_{m},m,\epsilon_{m}),T\big)$
and $\alpha_{m}((T,S))$ for the domino tableau pair obtained
by following the above shape-tracking procedure for $\alpha_{m}(T)$.

When $r=0$ or $1$,  the $G_r$ are precisely Garfinkle's algorithms;
for $r>1$ they are natural extensions to larger-rank tableaux. In
all cases, $G_r$ defines a bijection from $H_n$ to pairs of
same-shape tableaux in $SDT_r(n)$ \cite{vanleeuwen:rank}. These
generalizations of the Robinson-Schensted algorithm share a number of
properties with the original algorithm.  We state the following:

\begin{proposition}{\sc (\cite{vanleeuwen:rank}, (4.2))} \label{prop} $G_r(\sigma^{-1})=(S,T)$
whenever $G_r(\sigma)=(T,S)$. In particular, if $\sigma$ is an
involution, $G_r(\sigma)=(T,T)$ for some standard domino tableau
$T$.
\end{proposition}

\begin{example}
Consider the signed permutation $(2 \; -4 \; -3 \; 1)$.  It corresponds to the set $\sigma=\{(2,1,1),(4,2,-1),(3,3,-1),(1,4,1)\} \in H_4$.  If $r=2$, then successive insertion of elements of $\sigma$ into the empty tableau of rank zero yields
the following sequence of tableau pairs

$$
\raisebox{1ex}{$T(1)=$ \;}
\begin{tiny}
\begin{tableau}
:.{0}.{0}>2\\
:.{0}\\
\end{tableau}
\end{tiny}
\hspace{1in} \raisebox{1ex}{$S(1)=$ \;}
\begin{tiny}
\begin{tableau}
:.{0}.{0}>1\\
:.{0}\\
\end{tableau}
\end{tiny}
$$
$$
\raisebox{1ex}{$T(2)=$ \;}
\begin{tiny}
\begin{tableau}
:.{0}.{0}>2\\
:.{0}\\
:^4\\
\end{tableau}
\end{tiny}
\hspace{1in} \raisebox{1ex}{$S(2)=$ \;}
\begin{tiny}
\begin{tableau}
:.{0}.{0}>1\\
:.{0}\\
:^2\\
\end{tableau}
\end{tiny}
$$

$$
\raisebox{1ex}{$T(3)=$ \;}
\begin{tiny}
\begin{tableau}
:.{0}.{0}>2\\
:.{0}^4\\
:^3\\
\end{tableau}
\end{tiny}
\hspace{1in} \raisebox{1ex}{$S(3)=$ \;}
\begin{tiny}
\begin{tableau}
:.{0}.{0}>1\\
:.{0}^3\\
:^2\\
\end{tableau}
\end{tiny}
$$

$$
\raisebox{1ex}{$T(4)=$ \;}
\begin{tiny}
\begin{tableau}
:.{0}.{0}>1\\
:.{0}>2\\
:^3>4\\
\end{tableau}
\end{tiny}
\hspace{1in} \raisebox{1ex}{$S(4)=$ \;}
\begin{tiny}
\begin{tableau}
:.{0}.{0}>1\\
:.{0}^3^4\\
:^2\\
\end{tableau}
\end{tiny}
$$

\vspace{.1in}
\noindent
Consequently,  $G_2(\sigma)=(T(4),S(4))$.

\end{example}

\subsection{Cycles}
The notion of a cycle in a domino tableau appears in a number of
references. See for instance \cite{carre-leclerc},
\cite{vanleeuwen:edge}, or \cite{vanleeuwen:bijective}.  We now
review its definition.

\begin{definition} For a standard domino tableau $T$
of arbitrary rank $r$, we will call a square in position $(i,j)$
{\it fixed} when $i+j$ has the opposite parity as $r$, otherwise,
we will call it {\it variable}.
\end{definition}

It is possible to choose the sets of fixed and variable squares
differently, as in \cite{garfinkle1},(1.5.4); however, we refrain
from defining the more general possibilities as only this choice
will be necessary for our results.

If $T \in SDT_r(n)$, we will write $D(k,T)$ for the domino labeled
by the positive integer $k$ in $T$ viewed as a set of labeled
squares, and $supp \, D(k,T)$ will denote its underlying squares.
Write $label \, S_{i,j}$ for the label of the square $S_{i,j}$ in
$T$.  We extend this notion slightly by letting $label \, S_{i,j}
=0$ if either $i$ or $j$ is less than or equal to zero, and $label
\, S_{i,j} =\infty$ if $i$ and $j$ are positive but $S_{i,j}$ is not
a square in $T$.

\begin{definition}
Suppose that  $supp \, D(k,T)= \{S_{i,j},S_{i+1,j}\}$ or
$\{S_{i,j-1},S_{i,j}\}$ and the square  $S_{i,j}$ is fixed. Define
$D'(k)$ to be a domino labeled by the integer $k$ with $supp \,
D'(k,T)$ equal to
    \begin{enumerate}
        \item $\{S_{i,j}, S_{i-1,j}\}$     if $k< label \, S_{i-1,j+1}$
        \item $\{S_{i,j}, S_{i,j+1}\}$    if $k> label \, S_{i-1,j+1}$
    \end{enumerate}
Alternately, suppose that  $supp \, D(k,T)= \{S_{i,j},S_{i-1,j}\}$
or $\{S_{i,j+1},S_{i,j}\}$ and the square  $S_{i,j}$ is fixed.
Define $supp \, D'(k,T)$ to be
    \begin{enumerate}
        \item $\{S_{i,j},S_{i,j-1}\} $    if $k< label \, S_{i+1,j-1}$
        \item $\{S_{i,j},S_{i+1,j}\}$         if $k> label \, S_{i+1,j-1}$
    \end{enumerate}
\end{definition}

\begin{definition}
The cycle $c=c(k,T)$ through $k$ in a standard domino tableau $T$ is
a union of labels of $T$  defined by the condition that $l \in c$ if
either
    \begin{enumerate}
        \item $l=k$,
        \item $supp \, D(l,T) \cap supp \, D'(m,T) \neq
        \emptyset$ for some $m \in c$, or
        \item $supp \, D'(l,T) \cap supp \, D(m,T) \neq
        \emptyset$ for some $m \in c$.
    \end{enumerate}
\end{definition}

We will often identify the labels contained in the cycle with their underlying dominos.
For a standard domino tableau $T$ of rank $r$ and a cycle $c$ in  $T$,
we can define a domino tableau $MT(T,c)$ by replacing every domino
$D(l,T) \in c$ by the corresponding domino $D'(l,T)$.  That the resulting
tableau $MT(T,c)$ is standard follows from \cite{garfinkle1},
(1.5.27).  In general, the shape of $MT(T,c)$ will either equal the
shape of $T$, or one square will be removed (or added to the core)
and one will be added. The cycle $c$ is called {\it closed} in the former
case and {\it open} in the latter. For an open cycle $c$ of a tableau $T$, we will write
$S_b(c,T)$ and $S_f(c,T)$ for the squares that have been removed (or
added to the core) and added by moving through $c$; we will often abbreviate this notation to $S_b(c)$ and $S_f(c)$ when no confusion can result. Let $U$ be a set
of cycles in $T$. According to \cite{garfinkle1}, (1.5.29), the
order in which one moves through a set of cycles does not matter,
allowing us to unambiguously write $MT(T,U)$ for the tableau
obtained by moving-through all of the cycles in $U$.

We next define the set of cycles that it will be necessary to move
through to describe the relationship between $G_r$ and $G_{r+1}$.

For $T \in SDT_r(n)$, we will write $\delta=\delta(T)$ for the set
of squares $S_{i,j}$ that satisfy $i+j=r+2$.  These are the squares
with positive labels adjacent to the core of $T$. All are variable
in our choice of fixed and variable squares. In order to obtain a
domino tableau of rank $r+1$, it will be necessary to clear all of
the squares in $\delta$. Simply moving through $\Delta(T)$, the
cycles in $T$ that pass through $\delta$, will achieve this effect.
However, when applied to a pair of tableaux of the same shape, the
resulting pair of tableaux may not be of the same shape.  To this
effect,  we would like to define a minimal set of cycles in a pair
of domino tableaux that will ensure this. More precisely, for a pair
$(T,S)$, we would like to find sets of cycles $\gamma = (\gamma(T) ,
\gamma(S))$ in both $T$ and $S$ with $\Delta(T) \subset \gamma(T)$
and $\Delta(S) \subset \gamma(S)$ such that $MT(T, \gamma(T))$ and
$MT(S, \gamma(S))$ have the same shape.

The natural notion to consider is an extended cycle
(\cite{garfinkle2}, (2.3.1)), which we now reconstruct.

\begin{definition}
Consider $(T,S)$ a pair of same-shape domino tableaux, $k$ a label
of a domino in $T$, and $c$ the cycle in $T$ through $k$.  The extended cycle
$\tilde{c}$ of $k$ in $T$ relative to $S$ is a union of cycles in
$T$ which contains $c$.  Further, the union of two cycles $c_1 \cup
c_2$ lies in  $\tilde{c}$ if either is contained in $\tilde{c}$ and,
for some cycle $d$ in $S$,  $S_b(d)$ coincides with a square of
$c_1$ and $S_f(d)$ coincides with a square of $MT(T,c_2)$.  The
symmetric notion of an extended cycle in $S$ relative to $T$ is
defined in the natural way.
\end{definition}

Let
$\tilde{c}$ be an extended cycle in $T$ relative to $S$. According
to the definition, it is possible to write $\tilde{c} =c_1 \cup
\ldots \cup c_m$ and find cycles $d_1, \ldots, d_m$ in $S$ such that
$S_b(c_i)=S_b(d_i)$ for all $i$, $S_f(d_m)=S_f(c_1)$, and
$S_f(d_i)=S_f(c_{i+1})$ for $1\leq i < m$. The union $\tilde{d}=d_1
\cup \dots \cup d_m$ is an extended cycle in $S$ relative to $T$
called the {\it extended cycle corresponding to} $\tilde{c}$.
Symmetrically, $\tilde{c}$ is the extended cycle corresponding to
$\tilde{d}$.

It is now possible to define a moving through operation for a pair
of same-shape domino tableaux.  If we let $b$ be the ordered pair
$(\tilde{c},\tilde{d})$ of extended cycles in $(T,S)$ that
correspond to each other, then we define
$$MT((T,S), b)= (MT(T,\tilde{c}),MT(S,\tilde{d})).$$
As desired, this operation produces another pair of same-shape
domino tableaux (\cite{garfinkle2}, (2.3.1)).  If $B$ is a family of ordered pairs  of extended cycles that correspond to each other, then  we can unambiguously define $MT((T,S), B)$, the operation of moving through all of the pairs simultaneously.

\section{A Domino Tableau Correspondence}

From the definitions of the previous section, it is apparent that
moving through all of the extended cycles that pass thorough
$\delta(T)$ and $\delta(S)$ of a same-shape domino tableau pair $(T,S)$ will not only increase the rank of the resulting tableau pair by one, but the two tableaux will also be of the same shape.  What is perhaps surprising is that this map, which
merely evacuates $\delta$ in the simplest manner that will keep the
domino tableau pair of the same shape, describes the relationship
between the Robinson-Schensted maps $G_r$ and $G_{r+1}$.

\subsection{Main Theorem}

We first simplify our notation slightly.  Consider a pair of domino
tableaux $(T,S)$ of rank $r$ and define $\gamma(T)$ to be the set of
extended cycles in $T$ through $\delta(T)$ relative to $S$.
 Similarly, let $\gamma(S)$  be the set of extended cycles in $S$ through
$\delta(S)$ relative to $T$.  If we write $\gamma$ for the ordered
pair of sets of extended cycles $(\gamma(T),\gamma(S))$, then let
$$MMT((T,S))= MT((T,S),\gamma)$$
be the minimal moving through map that clears all of the squares in
$\delta(T)$ and $\delta(S)$.

\begin{theorem} \label{theorem:main} Consider an element $\sigma \in H_n$.  The Robinson-Schensted
  maps  $G_r$ and $G_{r+1}$ for rank $r$ and $r+1$ domino tableaux are related by
$$G_{r+1}(\sigma) = MMT(G_{r}(\sigma)).$$
\end{theorem}

The  proof is a direct consequence of the following lemma;  we show
that domino insertion commutes with moving through the set of
extended cycles which pass through the squares adjacent to the cores
of a domino tableau pair. We note that the lemma is not true when
more general sets of cycles are considered.

\begin{lemma} \label{lemma:main} Consider $\sigma \in H_n$.  Then

$$MMT\big(\alpha_{k+1}(G_r^k(\sigma))\big)=
\alpha_{k+1}\big(MMT(G_r^k(\sigma))\big)
$$
\end{lemma}

When $r=0$, the result is reminiscent of \cite{garfinkle2}, (2.3.2).
We follow a similar approach and redefine the scope of a number of
technical statements to cover the situations possible in the set of
rank $r$ standard domino tableaux when $r \geq 0$.

\begin{example}
Consider $\sigma=((2,1,-1),(1,2,1))$ in $H_2$.  If
$(T,S)=G_0(\sigma)$, then

$$
\raisebox{1ex}{$T=$ \;}
\begin{tiny}
\begin{tableau}
:>1\\
:>2\\
\end{tableau}
\end{tiny}
\hspace{1in} \raisebox{1ex}{$S=$ \;}
\begin{tiny}
\begin{tableau}
:^1^2\\
:;\\
\end{tableau}
\end{tiny}
$$
The cycles in $T$ are $c_1=\{1\}$ and $c_2=\{2\}$ and the cycles in
$S$ are $d_1=\{1\}$ and $d_2=\{2\}$.   Note that $\Delta(T)=c_1$ and
$\Delta(S)=d_1$.  However, $\gamma(T)=c_1 \cup c_2$ and
$\gamma(S)=d_1 \cup d_2$, so that $MMT(G_0(\sigma))$ is the pair of
tableaux

$$
\raisebox{2ex}{$T'=$ \;}
\begin{tiny}
\begin{tableau}
:.0>1\\
:^2\\
:;\\
\end{tableau}
\end{tiny}
\hspace{1in} \raisebox{2ex}{$S'=$ \;}
\begin{tiny}
\begin{tableau}
:.0>2\\
:^1\\
:;\\
\end{tableau}
\end{tiny}
$$
As stated in the theorem, $MMT(G_0(\sigma))\equiv(T',S')$ equals
$G_1(\sigma)$.
\end{example}

\subsection{Technical Lemmas}

It is possible to describe the open cycles in $T(k+1)$ in terms of
the open cycles in $T(k)$.  Garfinkle's \cite{garfinkle2}, (2.2.3)
describes this relationship when $r=0$. With only minor changes,
this result can be stated for arbitrary rank tableaux. We will
write $OC(T)$ for the set of open cycles in $T$. To be precise,
let us recall a definition:

\begin{definition}
If $T_1, T_2 \in SDT_r(n)$, and $U_1$ and $U_2$ are sets of open
cycles in $T_1$ and $T_2$, then a map $\mu:U_1 \rightarrow U_2$ is a
{\it cycle structure preserving bijection} if for every $c \in U_1$,
$S_b(\mu(c))=S_b(c)$ and $S_f(\mu(c))=S_f(c)$.
\end{definition}

In general, there is no cycle structure preserving bijection between
the open cycles in  $T(k+1)$ and those in $T(k)$.  However, their
relationship is only slightly more subtle.

\begin{definition}
A cycle $c \in OC(T(k+1))$ corresponds to a cycle $c' \in OC(T(k))$
if either $S_b(c')=S_b(c)$ or $S_f(c)=S_f(c')$.
\end{definition}

We will describe the open cycle correspondences and cycle structure
preserving bijections between $T(k+1)$ and $T(k)$.  The first lemma
is a generalized version of \cite{garfinkle2}, (2.2.3), extended by
the case here labeled as $2(a)(ii)$. Before stating it, let us
introduce notation that will be used throughout this section. We will
write $T$ for $T(k+1)$, $T'$ for $T(k)$, and $\overline{U}$ for the
tableau $U$ with its highest-labeled domino removed.  Let $P$ be the
squares in $T$ that are not in $T'$ and $\overline{P}$ be the
squares in $\overline{T}$ that are not in $\overline{T}'$.  If $e$
is the highest label in $T$, let $P_e'$ be the squares of $D(e,T')$,
and $P_e$ be the squares of $D(e,T)$.

\begin{lemma}\label{lemma:technical}
Consider $T(k)$ and $T(k+1) \in SDT_r(n)$.  Suppose $P$ is
horizontal and consists of the squares $\{S_{ij},S_{i,j+1}\}$.  When
$P$ is vertical instead, the obvious transpositions of the below
statements are true. The relationship of the open cycle structure of
$T(k)$ to the open cycle structure of $T(k+1)$ is described by the
following cases:
    \begin{enumerate}
        \item Suppose $S_{i,j+1}$ is variable.
            \begin{enumerate}
                \item First assume that $j>1$ and $S_{i+1,j-1}$  is not contained in the diagram underlying
                $T(k)$.  Let $c'$ be the open cycle  in $T(k)$
                with $S_b(c')=S_{i,j-1}$. Then there is an open cycle $c$ in  $T(k+1)$
                with $S_f(c)=S_f(c')$ and $S_b(c)=S_{i,j+1}$.
                Furthermore, there is a cycle structure preserving bijection
                between the remaining open cycles of $T(k)$ and
                $T(k+1)$.

                \item Otherwise, either $j=1$ or $S_{i+1,j-1}$ is contained in the diagram underlying
                $T(k)$. Then there are two possibilities.  Either
                    \begin{enumerate}[(i)]
                        \item there is an open cycle $c$ in  $T(k+1)$
                        with $S_b(c)=S_{i,j+1}$ and $S_f(c)=S_{i+1,j}$
                        and a cycle structure preserving bijection between $OC(T(k))$ and
                        $OC(T(k+1))\setminus \{c\}$, or

                        \item there is an open cycle $c'$ in
                        $T(k)$ and cycles $c_1, c_2$ in  $T(k+1)$
                        such that $S_f(c_1)=S_f(c')$,
                        $S_b(c_1)=S_{i,j+1}$, $S_f(c_2)=S_{i+1,j}$, and
                        $S_b(c_2)=S_b(c')$.  In this case, there is a cycle structure preserving bijection between $OC(T(k))\setminus \{c'\}$ and
                        $OC(T(k+1))\setminus \{c_1,c_2\}$.
                    \end{enumerate}
            \end{enumerate}

             \item Suppose $S_{i,j+1}$ is fixed.
                \begin{enumerate}
                    \item First assume that either $i=1$ or
                    $S_{i-1,j+2}$ is contained in the diagram underlying
                    $T(k+1)$.  There are two possibilities.  Either

                    \begin{enumerate}[(i)]
                        \item there is an open cycle $c'$ in $T(k)$ with
                        $S_f(c')=S_{ij}$ and an open cycle
                        $c$ in  $T(k+1)$ with $S_f(c)=S_{i,j+2}$ and  $S_b(c)=S_b(c')$;
                        in this case there is a cycle structure preserving bijection
                        between the remaining open cycles of $T(k)$ and
                        $T(k+1)$, or

                        \item $S_{ij} \in \delta(T(k))$,
                        there is a cycle $c$ in  $T(k+1)$ with
                        $S_b(c)=S_{i,j}$ and $S_f(c)=S_{i,j+2}$, and a cycle structure preserving bijection
                        between $OC(T(k))$ and $OC(T(k+1))\setminus
                        \{c\}$.
                    \end{enumerate}

                    \item Otherwise, both $i>1$ and $S_{i-1,j+2}$ is not contained in the diagram underlying
                    $T(k+1)$.
                      Then there is an integer $u >
                    \sigma_{k+1}$ such that the domino with label
                    $u$ forms a cycle $c'$ in $T(k)$ with
                    $S_f(c')=S_{ij}$ and $S_b(c')=S_{i-1,j+1}$.  In
                    this case,
                    there is a cycle structure preserving bijection between $OC(T(k))\setminus
                    \{c'\}$ and $OC(T(k+1))$.

                \end{enumerate}
    \end{enumerate}
\end{lemma}

To verify the above, it is necessary to understand how the cycle
structure of a domino tableau $U$ is related to the cycle structure
of $\overline{U}$. When $r=0$, this is described in
\cite{garfinkle2}, (2.2.4). Again for completeness, we state our
version for arbitrary rank tableaux in full, which differs in the
additional case $2(a)(ii)$. The proof of this lemma follows from an
easy, but tedious, inspection.

\begin{lemma}\label{lemma:2}
    Suppose that $T \in SDT_r(n)$, $e$ is the label of its highest domino $D$, and
    $\overline{T}$ is the domino tableau  with $D$ removed.
    Suppose $D$ occupies the squares
    $\{S_{ij},S_{i,j+1}\}$ in $T$.  Again,
    the obvious transpositions of the statements below are true for
    vertical $D$.
        \begin{enumerate}
        \item Suppose that $S_{i,j+1}$ is variable.
            \begin{enumerate}
                \item First assume that $j>1$ and $S_{i+1,j-1}$  is not contained in the diagram underlying
                $\overline{T}$. Let $\overline{c}$
                be the open cycle in $\overline{T}$ with
                $S_b(\overline{c}) = S_{i,j-1}$.  Then
                there is an open  cycle $c$ in  $T$ with
                $S_f(c)=S_f(\overline{c})$ and
                $S_b(c)=S_{i,j+1}$.
                Furthermore, there is a cycle structure preserving bijection
                between the remaining open cycles of $\overline{T}$ and
                $T$.

                \item Otherwise, either $j=1$ or $S_{i+1,j-1}$  is  contained in the diagram underlying
                $\overline{T}$.  Then $c=\{e\}$ is an open cycle in
                $T$ and there is a cycle structure preserving bijection
                between  $OC(\overline{T})$ and $OC(T)\setminus \{c\}$.
            \end{enumerate}

             \item Suppose that $S_{i,j+1}$ is fixed.
                \begin{enumerate}
                    \item First assume that  either $i=1$ or $S_{i-1,j+2}$ is contained in the diagram underlying
                    $T$.Then there are two possibilities.  Either

                    \begin{enumerate}
                        \item there exists an open cycle  $\overline{c}$ in $\overline{T}$
                        with $S_f(\overline{c})=S_{ij}$, and $c= \overline{c} \cup
                        \{e\}$ is an open cycle in $T$; in this case there is a cycle structure preserving bijection
                        between
                        $OC(\overline{T})\setminus \{\overline{c}\}$ and $OC(T)\setminus
                        \{c\}$, or

                        \item $S_{ij} \in \delta(T)$,
                        there is a cycle $c$ in  $T$ with
                        $S_b(c)=S_{ij}$ and $S_f(c)=S_{i,j+2}$, and a
                        cycle structure preserving bijection
                between $OC(\overline{T})$ and $OC(T)\setminus
                        \{c\}$.
                    \end{enumerate}

                    \item Otherwise, both  $i>1$ and  $S_{i-1,j+2}$ is not contained in the diagram underlying
                    $T$. Then either
                        \begin{enumerate}
                            \item there is a cycle $\overline{c}$ in
                            $\overline{T}$ with
                            $S_b(\overline{c})=S_{i-1,j+1}$ and
                            $S_f(\overline{c})=S_{ij}$,
                            $c= \overline{c} \cup \{e\}$ is a closed cycle
                            in $T$, and $OC(T)=OC(\overline{T})
                            \setminus\{\overline{c}\}$, or

                            \item there are two open cycles
                            $\overline{c}_1, \overline{c}_2$ in $\overline{T}$ such
                            that $S_b(\overline{c}_1)= S_{i-1,j+1}$,
                            $S_f(\overline{c}_2)=S_{ij}$, the set
                            $c=c_1 \cup c_2 \cup \{e\}$ is an open cycle
                            in $T$ and $OC(T)\setminus\{c\} =
                            OC(\overline{T})\setminus \{\overline{c}_1,\overline{c}_2\}$.
                        \end{enumerate}
                \end{enumerate}
    \end{enumerate}
\end{lemma}

\noindent Armed with this observation, we can now prove Lemma
\ref{lemma:technical}.
\begin{proof}

Lemma \ref{lemma:2} describes the relationships between the cycle
structures of  $\overline{T}(k)$ and $T(k)$, as well as
$\overline{T}(k+1)$ and $T(k+1).$  If we use induction on the size
of the tableaux, we can relate the cycle structures of
$\overline{T}(k)$ and $\overline{T}(k+1)$.  Together, this allows us
to describe the desired relationship between the cycle structures of
$T(k)$ and $T(k+1)$.

 If a pair of squares in a domino
tableau satisfy the hypotheses of a case of Lemma
\ref{lemma:technical} or Lemma \ref{lemma:2}, we will say that the
pair lies in the situation labeled by that case.  The proof of the
lemma divides into different cases described by the situations of
$\overline{P}$ and $P_e'$ and their relative positions.  When $r=0$,
this is exhaustively carried out in the proof of \cite{garfinkle2},
(2.2.3), which includes a description of the possibilities for
$\overline{P}$ and $P_e'$.  We will use the same labels for these
possibilities.  To verify the lemma for arbitrary rank tableaux, we
must check that the conclusions still hold in the cases originally
considered, as well as examine the new cases that arise for larger
rank tableaux. The former follows from a lengthy inspection of the
proof of \cite{garfinkle2}, (2.2.3).  We examine the new cases.

We have to consider situations where either $P, \overline{P}, P_e',$
or $P_e$ is in situation $2(a)(ii)$.  Most of the cases are
essentially trivial.  We treat two of them in detail; the rest
follow along similar lines.  The cases are labeled to mimic similar
cases considered in \cite{garfinkle2}, (2.2.3).

 Case K$'$. Here $\overline{P}=P_e'$ is in
situation $2(a)(ii)$. We have a cycle structure preserving bijection
between $OC(\overline{T})$ and $OC(T')$.  Note that $P_e=P$, and
they both must be in situation $2(a)(ii)$ or $1(b)$.  In both cases,
the desired relationship between $OC(T)$ and $OC(T')$ exists between
$OC(T)$ and $OC(\overline{T})$ by Lemma \ref{lemma:2}. Since we already
have a cycle structure preserving bijection between
$OC(\overline{T})$ and $OC(T')$, we are done.

Case L$'$. Here $\overline{P}$ is in situation $2(a)(ii)$ and
$P_e'=\{S_{i,j},S_{i+1,j}\}$, so that $P_e'$ is in situation
$2(a)(ii)$ as well.  If $D$ is the domino in $\overline{T}$ in
position $\overline{P}$ with label $f$, then we have a cycle
structure preserving bijection between $OC(\overline{T}) \setminus
\{f\}$ and $OC(\overline{T}') \setminus \{e\}$. Note that
$P=\{S_{i,j+1},S_{i+1,j+1}\}$ is in situation $1(b)$ of
\ref{lemma:technical} and $P_e=\{S_{i+1,j},S_{i+1,j+1}\}$ is in
situation $1(b)$ of Lemma \ref{lemma:2}. Because of the latter, we know
there is a cycle structure preserving bijection between
$OC(\overline{T})$ and $OC(T) \setminus \{e\}.$ From this, we can
construct a cycle structure preserving bijection between $OC(T')
\setminus \{c'\}$ and $OC(T) \setminus \{c_1, c_2\}$ where
$c'=\{e\}$ in $T'$, $c_1=\{e\}$ in $T$, and $c_2 =\{f\}$, as
required in the conclusion of $1(b)(ii)$.

\end{proof}

\begin{lemma}\label{lemma:gamma}
The set $\gamma(T(k+1))$ is the union of the open cycles that
correspond to cycles in $\gamma(T(k))$ and the cycles through
$\delta(T(k+1))$.
\end{lemma}

\begin{proof}
Let us write $\tilde{\gamma}(T')$ for the set of open cycles in $T$
that correspond to open cycles in $\gamma(T')$. We may take $k>1$,
otherwise this is trivial. First assume that $\sigma_{k+1} =e$. Then
$P_e =\{S_{1,s},S_{1,s+1}\}$ and could be in situations $1(a)$,
$1(b)$, $2(a)(i)$, or $2(a)(ii)$ of Lemma \ref{lemma:2}.  In the first and third
cases, let $c'$ be the cycle in $T'$ through the square $S_{1,s-1}$.
Then $c=c'\cup \{e\}$ is the open cycle in $T$ corresponding to
$c'$, $OC(T')\setminus \{c'\} = OC(T) \setminus \{c\}$, and $c \in
\gamma(T)$ iff $c' \in \gamma(T')$. Since $\Delta(T) \subset
\tilde{\gamma}(T')$, the result follows.  If $P_e$ is in situation
$1(b)$ of Lemma \ref{lemma:2}, then $OC(T')=OC(T)\setminus\{e\}$.
Since $\{k+1\}$ is a cycle in $S$, $\{e\}$ must be an extended cycle
implying that $\{e\} \notin \gamma(T)$.  Again, $\Delta(T) \subset
\tilde{\gamma}(T')$ and the result follows. If $P_e$ is in situation
$2(a)(ii)$ of Lemma \ref{lemma:2}, then $\{k+1\}$ is a cycle in $S$, $\{e\}$ must be an extended cycle and since $\{e\} \in \Delta(T)$, the result follows.

The rest of the proof is by induction on the size of the tableau. We
will assume that $\gamma(\overline{T})=\tilde{\gamma}(\overline{T}')
\cup \Delta(T)$.  We treat cases A-C and L from the proof of
\cite{garfinkle2}, (2.2.3) incorporating the additional
possibilities that arise in higher rank tableaux. Remaining cases
are handled along similar lines.

Case A. Suppose $\overline{P}$ is in situation $1(a)$ and
$\overline{P}=P_e'$.  Then $P= P_e$ and they both equal to the set
$\{S_{i+1,s},S_{i+1,s+1}\}$ for some $s$.  The squares of $P$ may be in
situations $1(a)$, $1(b)$, $2(a)(i)$, or $2(a)(ii)$ of
Lemma \ref{lemma:technical}.  In the first case, consider $c'$ as in Lemma
$\ref{lemma:technical}  (1(a)) $. The cycle $c'$ corresponds to
$c=c(e,T)$ since $S_f(c)=S_f(c')$. Examining the position of
$D(k+1,S)$, we find that the rest of the extended cycle
structure of $T$ is the same as in $T'$.  Hence if $c$ is any cycle
in $T$ that corresponds to a cycle $c'$ in $T'$, then $c \in
\gamma(T)$ iff $c' \in \gamma(T')$.  If $P$ lies in situation
$2(a)(ii)$, then $S_{i+1,s} \in \delta(T)$, $c(e,T)$ is a cycle
through $\delta(T)$ and lies in $\gamma(T)$.  Similar arguments work
for the remaining two cases.

Case B. Here $\overline{P}$ is in situation $1(a)$ and
$P_e'=\{S_{i+1,j-1},S_{i+2,j-1}\}$, implying that $P=\overline{P}$
and $P_e=P_e'$.  First consider the cycle $c=c(e,T)=\{e\}$.  Note
that $c$ corresponds to $c'=c(e,T')$ since $S_b(c)=S_b(c')$.  Let
$\overline{c}' = c' \setminus \{e\} \in OC(\overline{T}')$. Let
$f=label(S_{i,j+1},T)$ and note that the squares of $P$ form a
domino in $S$ with label $k+1$.  Then $S_b(k+1,S)=S_b(f,T)$ and
$S_f(k+1,S)=S_f(e,T)$, so that $e$ and $f$  are both in the same extended cycle of $T$ relative to $S$.  Hence $e \in \gamma(T)$
iff $f \in \gamma(T)$ iff $f \in \gamma(\overline{T})$ iff
$\overline{c}' \in \gamma(T')$ iff $c' \in \gamma(T')$, as desired.
For any open cycle $c$ not containing $e$ in $T$, the result follows
by induction.

Case C. Here $\overline{P}$ is in situation $1(a)$ and
$P_e'=\{S_{i+1,j-2},S_{i+1,j-1}\}$ is in situation $2(b)(i)$. Then
$P=\overline{P}$ and $P_e=P_e'$.  Let $c=c(e,T)$ and by the
conclusion of Lemma \ref{lemma:2} we find $S_f(c) = S_{i+1,j}$ and
$S_b(c)=S_{i,j+1}$. Note that $c$ corresponds to no open cycles in
$T'$. Since $S_f(c,T)=S_f(k+1,S)$ and $S_b(c,T)=S_b(k+1,S)$, the extended cycle of $e$ is just $c$.
Hence $c \in \gamma(T)$ iff $c$ passes through $\delta(T)$. For any
open cycle $c$ not containing $e$ in $T$, the result follows by
induction.

Case L. Consider $\overline{P}$  in situation $2(a)(ii)$ and
$P_e'=\{S_{i,j},S_{i+1,j}\}$, so that $P_e'$ is in situation
$2(a)(ii)$ as well.  We then have $P=\{S_{i,j+1},S_{i+1,j+1}\}$ and
$P_e=\{S_{i+1,j},S_{i+1,j+1}\}$.  First consider the cycle
$c=c(e,T)=\{e\}$. Note that $\overline{P}$ is a domino in $T$, say
with label $f$, and $P_e'$ is a domino in $S$, say with label $l$.
Then $S_b(c(l,S),S)=S_{i,j}=S_b(c(f,T),T)$ and
$S_f(c(l,S),S)=S_{i+2,j}=S_f(c(e,T),T)$.  Hence $\{e\}$ lies in the
extended cycle through $c(f,T)$.  Since $c(f,T) \in \Delta(T)$, we
must have $\{e\} \in \gamma(T)$.  If we let $c'=c(e,T')$, then
$S_f(c)=S_f(c')$, which means that $c$ corresponds to $c'$.  In
other words, $\{e\}$ lies in $\gamma(T)$ and $\tilde{\gamma}(T')
\cup \Delta(T)$.  Finally, consider any open cycle $c$ not
containing $e$ in $T$.  Then $c$ is also an open cycle in
$\overline{T}$, and the rest follows by induction. We omit the
argument when $\overline{P}$ and $P_e'$ are in situation $2(a)(i)$
instead.

\end{proof}

If we abuse notation and write $MMT(T)$ for $MT(T,\gamma(T))$, then
we can state the following version of Garfinkle's \cite{garfinkle2},
(2.2.9), which verifies Lemma \ref{lemma:main} for left tableaux.

\begin{lemma}\label{lemma:left}
Consider $\sigma \in H_n$  and write  $T(m)$  for the left tableau
of $G_r^m(\sigma)$. Then
$$\alpha_{k+1}\big(MMT(T(k))\big) = MMT\big(T(k+1)\big).$$
\end{lemma}

\begin{proof}
Using Lemma \ref{lemma:gamma},  we have to show that
$$\alpha_{k+1}\big(MMT(T(k))\big) = MT(T(k+1), \tilde{\gamma}(T(k))\cup\Delta(T(k+1))).$$
which is an adaptation of \cite{garfinkle2}, (2.2.9).  However, we
cannot adapt the proof of \cite{garfinkle2}, (2.2.9) verbatim, as it
uses induction on the number of open cycles in the extended cycle
defining the moving through operation.  In our situation, moving
through a set of cycles smaller than $\gamma(T(k))$ may leave us
with a domino tableau on which $\alpha$ is undefined. Nevertheless,
since only one pair $P$ of squares is added to $T(k)$ with domino
insertion, and moving through open cycles can be done independently,
we can essentially follow the original proof and examine the
relationship of $P$ with the cycles in $\gamma(T(k))$ individually.

The case when $\sigma_{k+1}=e$ is simple, and we assume that
$\sigma_{k+1} \neq e$.  We proceed by induction on $n$, noting that
the case $n=1$ corresponds to $\sigma_{k+1}=e$.  Following the
original proof of \cite{garfinkle2}, (2.2.9),  we show that each
domino in $\alpha_{k+1}\big(MMT(T')\big)$ lies in the same position
in $MMT(T)$.  For dominos with labels less than $e$, this will
follow by induction; for the domino with label $e$, it will follow
by inspection of each of the cases below.

Let $\overline{P}_1$ be the squares in $\alpha_{k+1}(\overline{T}')$
that are not in  $\overline{T}'$,  $\overline{P}_2$ be the squares
in $\alpha_{k+1}(MMT(\overline{T}'))$ that are not in
$MMT(\overline{T}').$ Write $T_1$ for $T$, $T_1'$ for $T'$, $T_2$
for $MMT(T)$, $T_2'$ for $MMT(T')$, and $T_3$ for
$\alpha_{k+1}(MMT(T'))$.  Hence we are verifying that $T_2=T_3$.

Case A.  Assume that $\overline{P}_1=P_e' = \{S_{ij},S_{i,j+1}\}$,
and $\overline{P}_1$ is in situation $1(a)$. Then $P_e = P=
\{S_{i+1,s},S_{i+1,s+1}\}$ for some $s$. Suppose first that
 $S_{i+1,s}$ is variable and that no cycle $c' \in \gamma(T_1')$ has
$S_f(c') = S_{i+1,s}$. If $S_{i+1,s} \in \delta$, then $\{e\} \in
\gamma(T_1)$ and $P_e(T_2)=\{S_{i+1,s+1},S_{i+1,s+2}\}=P_e(T_3)$.
When $S_{i+1,s} \notin \delta$, we have
$P_e(T_2)=P_e(T_1)=P_e(T_3)$. Suppose next that there is a cycle $c'
\in \gamma(T_1')$ with $S_f(c')=S_{i+1,s}$, then $e$ lies in a cycle
in $\gamma(T_1)$ and $P_e(T_2) = \{S_{i+1,s+1},S_{i+1,s+2}\} =
P_e(T_3)$. If $S_{i+1,s}$ is fixed, then
$P_e(T_2)=\{S_{i+1,s-1},S_{i+1,s} \}=P_e(T_3)$ if $S_{i+1,s-1}$ lies
in some cycle of $\gamma(T_1')$, and
$P_e(T_2)=\{S_{i+1,s},S_{i+1,s+1} \}=P_e(T_3)$ if it does not.

Case K$'$. Here $\overline{P}_1=P_e'$ are in situation $2(a)(ii)$.
Then $P_e=\{S_{i+1,j},S_{i+1,j+1}\}$. Note that $c=\{e\}$ is a cycle
in $T_1$  and $d=\{k+1\}$ is a cycle in $S(k+1)$ with
$S_f(c,T_1)=S_{i+2,j}=S_f(d,S(k+1))$ and $S_b(c,T_1)=S_{i+1,j+1}=S_b(d,S(k+1))$. Hence $c=\{e\}$ is an extended cycle not contained in $\gamma(T_1)$ and consequently
$D(e,T_2)=D(e,T_1)=\{S_{i+1,j},S_{i+1,j+1}\}.$ Now note that
$\overline{P}_2=P(e,T_2')$ and by a similar argument, we obtain
$D(e,T_3)=\{S_{i+1,j},S_{i+1,j+1}\},$ as desired.

Case L$'$. Here $\overline{P}_1$ is in situation $2(a)(ii)$ and
$P_e'=\{S_{ij},S_{i+1,j}\}$, so it is in situation in $2(a)(ii)$ as
well.  Then $P_e=\{S_{i+1,j},S_{i+1,j+1}\}$ and
$P=\{S_{i,j+1},S_{i+1,j+1}\}$. Note that $c'=\{e\}$ is a cycle in
$T_1'$ with $S_f(c)=S_{i+2,j}$ and that the squares $P_e'$ form a
domino in $S(k)$, say with label $f$. Let $d = c(f, S(k+1))$ and
note $d \in \gamma(S(k+1))$. Furthermore,  $S_f(d)=S_{i+2,j}$
implying that $c(e,T_1) \in \gamma(T_1)$, and
$D(e,T_2)=\{S_{i+1,j},S_{i+2,j}\}$.  Now observe that
$D(e,T_2')=\{S_{i+1,j},S_{i+2,j}\}$ and
$\overline{P}_2=\{S_{i,j+1},S_{i,j+2}\}$.  This means $D(e,T_3) =
\{S_{i+1,j},S_{i+2,j}\}$, and $D(e,T_2)=D(e,T_3)$, as desired.

\end{proof}
\subsection{Domino Insertion and Moving Through}

Armed with the technical results of the previous section, we can now address the main lemma of the paper.  We prove Lemma \ref{lemma:main}, verifying that domino insertion on tableau pairs commutes with the minimal moving through map.   Write $(T_1',S_1')= (T(k),S(k))$,
$(T_1,S_1)=(T(k+1),S(k+1))$, $(T_2',S_2')=MMT(T_1',S_1')$,
$(T_2,S_2)= \alpha_{k+1}((T_2',S_2'))$, and $(T_3,S_3)=
MMT(T_1,S_1)$. Expressed in this notation, we would like to prove
that $(T_2,S_2)=(T_3,S_3)$.  Lemma \ref{lemma:left} says that
$T_2=T_3$, and it remains to show that $S_2=S_3$.

\begin{proof}
Write $P_1$ for the squares in $T_1$ that are not in $T_1'$ and
$P_2$ for the squares in $T_2$ that are not in $T_2'$.  Note that
$P_1$ forms a domino in $S_1$  and $P_2$ forms a domino in $S_2$,
both with label $k+1$.  Assume that $P_1=\{S_{i,j},S_{i,j+1}\}$.  We
will examine the cases when $P_1$ is in situations $1(a)$,
$1(b)(ii)$, and $2(b)$.  The others follow along similar lines.

So suppose that $P_1$ is in situation $1(a)$ of Lemma
\ref{lemma:technical} and $c$ is the open cycle with
$S_b(c)=S_{i,j+1}$ described therein. Then $D(k+1,S_1)$ is in
situation $1(a)$ of Lemma \ref{lemma:2} and there is an open cycle
$\overline{d}$ in $S_1'$ with $S_b(\overline{d})=S_{i,j-1}$ such
that $d=\overline{d}\cup \{k+1\}$ is an open cycle in $S_1$.  Note
that $c \in \gamma(T_1)$ iff $d \in \gamma(S_1)$.   If $c \in
\gamma(T_1)$, then by Lemma \ref{lemma:left}, $P_2=\{S_{i,j-1},S_{ij}\}$,
which implies that $D(k+1,S_3)=D(k+1,S_2)$.  Since the rest of the
cycle structure in $T_1$ remains the same as in $T_1'$, the rest of
the cycles in $\gamma(S_1)$ are the same as in $\gamma(S_1')$ and
consequently, $S_2=S_3$. If $c \notin \gamma(T_1)$, the result is
clear.

If $P_1$ is in situation $1(b)(ii)$ of Lemma \ref{lemma:technical},
then $D(k+1,S_1)$ is in situation $1(b)$ of Lemma \ref{lemma:2}. Let
$c',c_1,$ and $c_2$ be as described in Lemma \ref{lemma:technical}
$1(b)(ii)$ and let $d=c(k+1,S_1)$.  Since $S_b(c_1)=S_b(d)$ and
$S_f(c_2)=S_f(d)$, $c_1$ and $c_2$ lie in the same extended cycle
relative to $d$, so $c_1, c_2 \in \gamma(T_1)$ iff $d \in
\gamma(S_1)$.  If $c_1, c_2 \in \gamma(T_1),$ then by Lemma
\ref{lemma:left}, $P_2=\{S_{ij},S_{i+1,j}\}$.  Since $d \in
\gamma(S_1)$, this means $D(k+1,S_2)=D(k+1,S_3)$.  Since the rest of
the cycle structure in $T_1$ remains the same as in $T_1'$, the rest
of the cycles in $\gamma(S_1)$ are the same as in $\gamma(S_1')$ and
we can conclude that $S_2=S_3$. If $c_1, c_2 \notin \gamma(T_1),$
the conclusion is the same.

The most troublesome case is when $P_1$ is in situation $2(b)$ of
Lemma \ref{lemma:technical}. Then $D(k+1,S_1)$ is either in
situation $2(b)(i)$ or $2(b)(ii)$ of Lemma \ref{lemma:2}.  So
suppose first that $D(k+1,S_1)$ is in situation $2(b)(i)$. Let
$\overline{d}$ be the cycle in $S_1'$ with
$S_f(\overline{d})=S_{i,j}$ and $S_b(\overline{d})=S_{i-1,j+1}$.
Then $d=\overline{d} \cup \{k+1\}$ is a closed cycle in $S_1$ and
consequently does not lie in $\gamma(S_1)$.   Let $c'$ be the cycle
in $T_1'$ with $S_f(c')=S_{ij}$ and $S_b(c')=S_{i-1,j+1}$. Then $c'$
is the entire extended cycle in $T_1'$ that corresponds to
$\overline{d}$ in $S_1'$; in particular, this means that $c' \notin
\gamma(T_1')$ and $\overline{d} \notin \gamma(S_1')$.  Consequently,
$S_2=S_3$.

Finally, consider $D(k+1,S_1)$ in situation $2(b)(ii).$ Let $d_1$
and $d_2$ be the cycles in $S_1'$ with $S_b(d_1)=S_{i-1,j+1}$ and
$S_f(d_2)=S_{ij}$.  Then $d_1 \cup d_2 \cup \{k+1\}$ is an open
cycle in $S_1$.  Let $c'$ be as in Lemma \ref{lemma:technical} 2(b)
and note that $c' \in \gamma(T_1')$ iff $d_1,d_2 \in \gamma(S_1')$.
If $c' \in \gamma(T_1')$, then $P_2=\{S_{i-1,j+1},S_{i,j+1}\}$ by
Lemma \ref{lemma:left} and we again conclude that $S_2=S_3$.   If
$c' \notin \gamma(T_1')$, the result is clear.

\end{proof}

\subsection{Restriction to Involutions}

 We follow van Leeuwen in the next  definition, which constructs a map
 between domino tableaux of unequal rank
\cite{vanleeuwen:rank}.

    \begin{definition} Let $r$ and $r'$ be non-negative integers and
    suppose that $T \in SDT_r(n)$.  We define the map $t_{r,r'}:SDT_r(n) \rightarrow
    SDT_{r'}(n)$  by setting $t_{r,r'}(T) = T'$ whenever
    $G_r^{-1}(T,T) = G_{r'}^{-1}(T',T')$.
    \end{definition}

Armed with Theorem \ref{theorem:main}, the maps $t_{r,r+1}$ take a
particularly simple form.  The domino tableau $t_{r,r+1}(T)$ in
$SDT_{r+1}(n)$ is simply the image of $T$ after all the cycles in
$\Delta(T)$ have been moved through.

\begin{corollary}
$t_{r,r+1}(T) = MT(T,\Delta(T))$
\end{corollary}
\begin{proof}
If $\sigma$ is an involution and $G_r(\sigma) =(T,S)$, then $S$ must
equal $T$.  The definition of extended cycles  implies that every
extended cycle in $T$ relative to $S$ consists of a unique cycle. In
our setting, this implies $\gamma = (\Delta(T),\Delta(T))$.   Using
Theorem \ref{theorem:main} and the definition of moving through
extended cycles, we now have that
$$\big(t_{r,r+1}(T),t_{r,r+1}(T)\big) = MT((T,T),\gamma)=\big(MT(T,\Delta(T)),MT(T,\Delta(T))\big),$$ as
desired.
\end{proof}

\end{document}